\documentclass[11pt,a4paper]{article}
\usepackage{epsfig}
\newtheorem{theorem}{Theorem}
\def\bkR{{\rm I\kern-.17em R}}
\begin{document}
\thispagestyle{empty}
%
%  Title
%
\title{\bf Poisson integrators for
Volterra lattice equations}
\author{
T. Ergen\c{c},
     B. Karas\"{o}zen\thanks{e-mail:
    bulent@metu.edu.tr}\\
          Department of Mathematics and Institute of Applied
     Mathematics\\
         Middle East Technical University, 06531 Ankara-Turkey}
%\\
%    {\small    $^{a}$ }
% Department of Mathematics,
%         Middle East Technical University, 06531 Ankara-Turkey}}
\date{}
\maketitle
\begin{abstract}
The Volterra lattice equations are completely integrable and possess
bi-Hamiltonian structure. They are integrated using partitioned
Lobatto IIIA-B methods which preserve the Poisson structure.
Modified equations are derived for the symplectic Euler and second
order Lobatto IIIA-B method. Numerical results confirm preservation
of the corresponding Hamiltonians, Casimirs, quadratic and cubic
integrals in the long-term with different orders of accuracy.
\end{abstract}

\medskip

%{  AMS classification:} 65L, 65M, 70H\\

\section{Introduction}The preservation of qualitative properties of discretized models of
continuous systems became more and more important in recent years.
The development of symplectic integrators for Hamiltonian systems,
construction of geometric integrators which preserve symmetries,
reversing symmetries and phase space volume of the underlying
differential equations are some examples. For Hamiltonian systems in
non-canonical form with a non-linear structure matrix, i.e Poisson
systems, there does not exist any general structure preserving
integrator similar to the symplectic methods for canonical
Hamiltonian systems (Sec. VII.2 \cite{hairer02gi},
\cite{suris97ano}) . There are some well developed methods based on
generating functions and Hamiltonian splitting for the Lie-Poisson
systems, i.e. non-canonical Hamiltonian systems with a linear
structure matrix (Sec. VII.2.6 \cite{hairer02gi}). Poisson systems
arise especially as Hamiltonian pde's like Korteweg de Vries (KdV)
and nonlinear Schr\"{o}dinger equation with infinitely many
integrals. Semi-discretization of these pde's in space by preserving
the integrals and Hamiltonian structure results in integrable
lattice equations; ode's with a Poisson structure and a finite
number of integrals.

After introducing the Volterra lattice equations in the next
section, we discuss in Section 3 symplectic integrators like the
implicit mid-point and the symplectic Euler method applied to the
well-known two dimensional Lotka-Volterra equation. Preservation of
the Poisson structure of the Volterra lattice equation by symplectic
Euler and Lobatto IIIA-B methods is proved and the corresponding
modified equations are derived. The numerical results show the
preservation of the conserved quantities by these methods.

\section{Volterra lattice equation}
The $m$-dimensional Volterra lattice
\begin{equation} \dot{y}_i = y_i(y_{i+1} - y_{i-1}), \quad
i=1,\ldots,m \label{VL}
\end{equation}
for even $m$ and with periodic boundary conditions $y_{m+i} = y_i,$
$i=0,1,\ldots$ and with $y_i > 0$ was studied first in
\cite{kac75oes} as an integrable system. It was shown that the
Volterra lattice equation represents an integrable discretization of
the KdV equation \cite{kac75oes, suris99idf} and of inviscid
Burger's equation \cite{kupershmidt97isd}. Besides these, the
Volterra lattice equation describes many phenomena such as the
vibrations of particles on lattices (Liouville model on the
lattice), waves in plasmas and the evolution of populations
 in a hierarchical system of competing species
\cite{cronstrom95mso, faddeev86lmo, suris99idf}.

The system (1) possesses  a bi-Hamiltonian structure
\cite{cronstrom95mso, suris99idf}
\begin{equation}
\dot{y} =  J_0(y)
\nabla H_1
 = J_1(y)
\nabla H_0
\end{equation}
with respect to the quadratic and the cubic Poisson brackets
\begin{equation}
\{ y_i,y_{i+1} \}_0 = y_iy_{i+1},
\end{equation}
\begin{equation}
\{ y_i,y_{i+1} \}_1 =  y_iy_{i+1} (y_i + y_{i+1}),\quad
\{ y_i,y_{i+2} \}_1 =  y_iy_{i+1} y_{i+2}.
\end{equation}
The corresponding Hamiltonians are
\begin{equation} H_1=
\sum_{i=1}^m y_i,\qquad H_0 = \frac{1}{2}\sum_{i=1}^m \log (y_i).
\label{HA}
\end{equation}
and the structure matrices $J_0(y)$ corresponding to the quadratic
Poisson bracket has the form
$$ J_0(y) = \left (
\begin{array}{ccccc}
0 & y_1y_2 &\ldots & \ldots &   -y_1y_m\\
\vdots & \ddots & \ddots & \ddots  & 0\\
\ldots & -y_{i-1}y_i  & 0 & y_iy_{i+1} & \ldots \\
\vdots & \vdots & \ldots &  \vdots  &  \vdots\\
y_1y_m & \ldots & \ldots  &  -y_{m-1}y_m & 0
\end{array} \right ).
$$
%$$
%J_1(y)  =
%$$
%\begin{small}
%$$
%\left (\begin{array}{ccccccc} 0 & y_1y_2(y_1+y_2) & y_1y_2y_3 &
%\ldots &  \ldots & \ldots
%-y_1y_{m-1}y_m & -y_1y_m(y_1+y_m) \\
%y_1y_2(y_1+y_2) &  -y_2y_3(y_2+y_3)
%& 0 & \ldots & \ldots & \ldots  -y_1y_2y_m \\
%\vdots & \ddots & \ddots & \ddots & \ddots & \ddots & \vdots \\
%\vdots & -y_{i-2}y_{i-1}y_i & -y_{i-1}y_i(y_{i-1}+y_i) & 0
%& y_iy_{i+1}(y_i+y_{i+1}) & y_iy_{i+1}y_{i+2} & \ldots\\
%-y_1y_{m-1}y_m & \ldots & 0 & y_{m-3}y_{m-2}y_{m-1} &
%-y_{m-2}y_{m-1}(y_{m-2}+y_{m-1})
%&  y_{m-1}y_m(y_{m-1}+y_m)\\
%\vdots & \ddots & \ddots & \ddots & \ddots & \ddots & \vdots  \\
%y_1y_m(y_1+y_m) & y_1y_2y_m & \ldots & \ldots  &
%-y_{m-2}y_{m-1}y_m & -y_{m-1}y_m(y_{m-1}y_m) & 0
%\end{array} \right )
%$$
%\end{small}

If the Poisson brackets are compatible, i.e. the sum of $\{ \cdot,
\cdot\}_0 + \{\cdot,\cdot\}_1$ is again a Poisson bracket, then the
bi-Hamiltonian system (1) has a finite number of functionally
independent first integrals $I_i,i=1,\ldots, m$ in involution i.e.
$\{ I_i, I_j\}_k=0$ for $k=0,1$ and $i\not = j$, i.e.. the Volterra
lattice is completely integrable with respect to both brackets.

The Hamiltonian $H_0$ is a Casimir with respect to the Poisson
bracket $\{\cdot,\cdot\}_0,$ i.e.$ \{H_0,F\}_0=0$ for any function
$F(y)$.

The Volterra lattice (1) represents an integrable discretization of
KdV equation  \cite{suris99idf}
\begin{equation}
\frac{\partial u}{\partial \tau} +6 u\frac{\partial u}{\partial
\xi} \frac{\partial^3 u}{\partial \xi^3}
 = 0, \label{KDV}
\end{equation}
which also possesses a bi-Hamiltonian structure  \cite{olver93aol}
and has infinitely many integrals. The first three integrals of
(\ref{KDV}) are \cite{frutos97aac, goktas98coc}
$$
{\cal I}_1(u)  =  \int_{-a}^{a} u d\xi,\quad {\cal I}_2(u)  =
\frac{1}{2}\int_{-a}^{a} u^2 d\xi,\quad {\cal I}_3(u)  =
\int_{-a}^{a} \left (\frac{1}{2} \left (\frac{\partial u}{\partial
\xi}\right) ^2- \frac{1}{6}u^3 \right ) d\xi.
$$
The corresponding conserved quantities of the Volterra lattice (1)
are the quadratic and cubic integrals
 \cite{goktas98coc}
\begin{equation}
I_q  =  \sum_{i=1}^m \frac{1}{2} y_i^2 + y_i y_{i+1},\quad I_c  =
\sum_{i=1}^m  \frac{1}{3} y_i^3 + y_i y_{i+1} \label{INT} (y_i +
y_{i+1} +y_{i+2}).
\end{equation}
In discretized form, the mass conservation ${\cal I}_1$
corresponds to the Hamiltonian $H_1$, the momentum and energy
integrals ${\cal I}_2,\;{\cal I}_3$ of the KdV equation (5)
correspond to the first integrals $I_q$ and $I_c$ of the Volterra
lattice (1) respectively.

The Volterra lattice is time-reversible and is closely connected to
the Toda lattice, which is also completely integrable and has a
tri-Hamiltonian structure \cite{damianou02ftl}, \cite{suris99idf}.
The transformation of variables $ a_i = y_{2i}y_{2i-1},\; b_i
=y_{2i-1} - y_{2i-2} $ in (1) gives the Toda lattice
$$
\dot{a}_i = a_i(b_{i+1} - b_i), \quad \dot{b}_i = a_i -
a_{i-1},\quad i=1,\ldots,m/2
$$
with periodic boundary conditions
$$
a_0 = a_{m/2}, \qquad b_{m/2+1} = b_1.
$$
\section{Poisson integrators}
Geometric integrators for the Poisson systems
\begin{equation}
\dot{y} = J(y) \nabla H \end{equation} with a skew-symmetric
non-constant structure matrix $J(y)$ were studied recently in
several papers. For recent surveys on Poisson integrators see (Sec.
VII.2, \cite{hairer02gi} and \cite{kar04pi}). For related material
on Poisson systems see (Sec. VII.2, \cite{hairer02gi}) , (Ch. 10,
\cite{marsden02ims}) and (Ch. 6 \& 7
 ,\cite{olver93aol}).

The Poisson bracket $\{F,G\}$ for two smooth functions $F(y)$ and
$G(y)$ is defined by
\begin{equation}
\{F,G\}(y) = \nabla F(y)^T J(y) \nabla G(y) \end{equation} which is
bilinear, skew-symmetric $\{G,F\} =-\{F,G\}$ and satisfies the
Leibniz' rule $(\{E\cdot F,G\} = E\cdot\{F,G\}+ F\cdot \{E,G\})$
 as well as the Jacobi identity $(\{E,\{F,G\}\} + \{F,\{G,E\}\} + \{G,\{E,F\}\}
 )$.

The structure matrix $J(y)$ of Poisson systems does not need to be
invertible as for canonical Hamiltonian systems with $J^{-1}$. All
odd-dimensional skew-symmetric structure matrices $J(y)$ are
singular. The structure matrices $J_0(y)$ and $J_1(y)$  of the
periodic Volterra lattice (1) are also singular. Such systems are
called degenerate Poisson systems.

The functions $C_i(y)$ satisfying $ \{C_i,H\}=0$ are called Casimirs
or distinguished functions, which are first integrals whatever
$H(y)$ is.

There are two characteristics for the flow $\varphi_t(y)$ of the
Poisson system (8): \begin{itemize}
\item the flow $\varphi_t(y)$ of the
differential equation (8) is a Poisson map, i.e.
$$
\dot{\varphi}_t(y) J(y) \dot{\varphi}_t(y)^T = J(\varphi_t(y))
$$
where $\dot{\varphi}_t(y)$ denotes the Jacobian of $\varphi$,
\item and it respects the Casimirs of $J(y)$, i.e.
$C_i(\varphi_t(y))= Const.$
\end{itemize}

A numerical method $y_{n+1} = \phi_h(y_n)$ is called a Poisson
integrator for the structure matrix $J(y)$ if the transformation
$y_n\rightarrow y_{n+1}$ is a Poisson map that respects the
Casimirs. The Casimirs should be preserved by the Poisson
integrator. But the Casimir functions can be arbitrary, in case of
the Volterra lattice they are in logarithmic form, therefore their
conservation depends on the special structure of the problem.

Because each Poisson system is distinguished by the structure matrix
$J(y)$, a method will be Poisson integrator only for a specific
class of structure matrices. Therefore symplectic methods used for
Hamiltonian systems can not be directly applied to Poisson systems.
But some of the symplectic integrators can preserve certain Poisson
structures. An example of this kind is the Poisson system resulting
from the Ablowitz-Ladik integrable discretization of the nonlinear
Schr\"odinger equation \cite{suris97ano} which is preserved by the
symplectic Euler method.

The two-dimensional Lotka-Volterra equation was studied by several
authors. The symplectic Euler and St\"ormer-Verlet methods preserve
the Poisson structure of (10) whereas the implicit mid-point rule
does not (see pp. 238 \cite{hairer02gi}). There are also some
non-standard methods(\cite{mickens03anf}, \cite{sanzserna94aus})
which preserve the Poisson structure of the two-dimensional
Lotka-Volterra equation.

In the following we will apply the symplectic Euler method and
Lobatto IIIA-B methods to the splitted form of the Volterra lattice
equation (1). Both  belong to the class of splitting methods which
have been successfully used as geometric integrators in recent
years. For a survey of the splitting methods see \cite{hairer03gni}
and \cite{mclachlan02sm}.

For the application of partitioned Runge-Kutta methods of Lobatto
type, we split the equation (1) into two parts
\begin{equation}
\dot{u}_i  = u_i(v_i-v_{i-1}),\quad
\dot{v}_i  =  v_i(u_{i+1}-u_i),\quad i=1,\ldots,m/2 \label{LVS}
\end{equation}
by grouping the variables into odd $u_i = y_{2i -1}$ and even $v_i =
y_{2i}$ parts. Equation (\ref{LVS}) is  bi-Hamiltonian with the
quadratic Poisson bracket and Hamiltonian like the Volterra lattice
\cite{suris99idf}:
\begin{equation}
\{u_i,v_i\}_0 = u_iv_i  \qquad  \{v_i,u_{i+1}\}_0 = v_iu_{i+1}.
\label{PB1}
\end{equation}
%\begin{eqnarray}
%\{u_i,u_{i+1}\}_1 & = &  u_iv_iu_{i+1},\; \{v_i,v_{i+1}\}_1 = v_iv_{i+1}u_{i+1}, \nonumber \\
%\{u_i,v_i\}_1  & = &  u_iv_i(u_i+v_i),\; \{v_i,u_{i+1}\}_1 = v_iu_{i+1}(v_i+u_{i+1})\label{PB2}
%\end{eqnarray}
The Hamiltonians and the first integrals \label{INT} can be
written in the new variables:
\begin{equation}
H_1(u) = \sum_{i=1}^{m/2} u_i + v_i, \quad H_0(u) =
\frac{1}{2}\sum_{i=1}^{m/2} \log (u_i)+ \log (v_i),
\end{equation}
\begin{equation}
 I_q  =  \sum_{i=1}^{m/2} \frac{1}{2} u_i^2 + u_i v_i,\quad I_c  =
\sum_{i=1}^{m/2}  \frac{1}{3} u_i^3 + u_i v_i  (u_i + v_i
+u_{i+1}).
\end{equation}

Symplectic and time-reversible partitioned Runge-Kutta methods like
the Lobatto IIIA-B methods \cite{hairer02gi} can be easily applied
to Volterra lattice equations in the partitioned form.  The
symplectic Euler method which consists of a combination of explicit
and implicit Euler is a first order  Lobatto IIIA-B method. For the
partitioned  Volterra equations (\ref{LVS}) the symplectic Euler
method becomes
\begin{equation}
u_i^{n+1}  =  u_i^n + h u_i^n(v_i^{n+1}-v_{i-1}^{n+1})  \qquad
v_i^{n+1}   =  v_i^n + h v_i^{n+1}(u_{i+1}^n-u_i^n). \label{SE}
\end{equation}

The second order Lobatto IIIA-B method for the partitioned
Volterra lattice (\ref{LVS}) gives
\begin{eqnarray*}
k^1_i   =  u_i^n \left(v_i^n + \frac{h}{2}l_i^1 -(v_{i-1}^n +\frac{h}{2}l_{i-1}^1) \right)
& & l_i^1 =  (v_i^n +\frac{h}{2}l_i^1)(u_{i+1}^n-u_i^n)\\
k^2_i  =  \left(u_i^n + \frac{h}{2}(k_i^1 +k_i^2)\right) \left(v_i^n + \frac{h}{2}l_i^1 -(v_{i-1}^n +\frac{h}{2}l_{i-1}^1)\right)
& & l_i^2 =  (v_i^n +\frac{h}{2}l_i^1)\left(u_{i+1}^n + \frac{h}{2}k_{i+1}^1 -(u_i^n +\frac{h}{2}k_i^1) \right)\\
u_i^{n+1} = u_i^n +  \frac{h}{2} (k_i^1 + k_i^2) & & v_i^{n+1} =
v_i^n +  \frac{h}{2} (l_i^1 + l_i^2).
\end{eqnarray*}
The internal stage vectors $k^2$ and $l^1$ are computed by solving a
system of linear equations, whereas the vectors $k^1$ and $l^2$ are
obtained explicitly. The second order Lobatto IIIA-B method is known
as St\"ormer-Verlet method for separable Hamiltonian systems.

In order to show the preservation of the Poisson structure with the
quadratic brackets (\ref{PB1}) we consider the corresponding
two-form formulation \cite{olver93aol}
\begin{equation}
\sum_{k=0}^{m/2} \frac{du_i\wedge dv_i}{u_iv_i} + \sum_{k=0}^{m/2}
\frac{du_i\wedge dv_{i+1}}{u_{i+1}v_i}. \label{PS1A}
\end{equation}
\begin{theorem}
The Poisson structure of the Volterra lattice with the quadratic
brackets (11) is preserved by the symplectic Euler method.
\end{theorem}
{\bf Proof:\/} In order to show the equality $$ \sum_{k=0}^{m/2}
\frac{du_i^{n+1}\wedge dv_i^{n+1}}{u_i^{n+1}v_i^{n+1}} +
\sum_{k=0}^{m/2} \frac{du_i^{n+1}\wedge
dv_{i+1}^{n+1}}{u_{i+1}^{n+1}v_i^{n+1}} = \sum_{k=0}^{m/2}
\frac{du_i^{n}\wedge dv_i^{n}}{u_i^{n}v_i^{n}} + \sum_{k=0}^{m/2}
\frac{du_i^{n}\wedge dv_{i+1}^{n}}{u_{i+1}^{n}v_i^{n}},
$$

we differentiate (11) to get
\begin{eqnarray}
du_i^{n+1} &  = &  du_i^n + h(v_i^{n+1} - v_{i-1}^{n+1})du_i^n + hu_i^ndv_i^{n+1} - hu_i^ndv_{i-1}^{n+1} \label{EU1},\\
dv_i^{n+1} &  = &  dv_i^n + h(u_{i+1}^n - u_i^n)dv_i^{n+1} + hv_i^{n+1}du_{i+1}^n - hv_i^{n+1}du_i^n \label{EV1}.
\end{eqnarray}

Using the equations (\ref{EU1}) and (\ref{EV1}) successively we
obtain
\begin{eqnarray*}
du_i^{n+1} \wedge dv_i^{n+1} & = & (1 +h(v_i^{n+1}-v_{i-1}^{n+1}))du_i^n\wedge dv_i^{n+1} - h u_i^ndv_{i-1}^{n+1}\wedge dv_i^{n+1}, \\
dv_i^{n+1} \wedge du_{i+1}^{n+1} & = & (1 +h(v_{i+1}^{n+1}-v_i^{n+1}))dv_i^{n+1}\wedge du_{i+1}^n + h u_{i+1}^ndv_i^{n+1}\wedge dv_{i+1}^{n+1}.
\end{eqnarray*}
The equations  (14) can be written in the equivalent form:
$$
1 + h (v_i^{n+1} - v_{i-1}^{n+1})  =  \frac{u_i^{n+1}}{u_i^n},\qquad
1 - h (u_{i+1}^n - u_i^n)  = \frac{v_i^n}{v_i^{n+1}}.
$$
Combining all these into the two-forms in (\ref{PS1A}) we obtain
\begin{eqnarray*}
du_i^{n+1} \wedge dv_i^{n+1} & = & \frac{u_i^{n+1} v_i^{n+1}}{u_i^nv_i^n}\left (du_i^n \wedge dv_i^n + h v_i^{n+1}du_i^n\wedge du_{i+1}^n\right )- hu_i^ndv_{i-1}^{n+1} \wedge dv_i^{n+1},\\
dv_i^{n+1} \wedge du_{i+1}^{n+1} & = & \frac{u_{i+1}^{n+1}v_i^{n+1}}{u_{i+1}^nv_i^n}\left (dv_i^n \wedge du_{i+1}^n - h v_i^{n+1}du_i^n\wedge du_{i+1}^n\right )+ hu_{i+1}^ndv_i^{n+1} \wedge dv_{i+1}^{n+1}.
\end{eqnarray*}
Because  $du_i^n \wedge du_{i+1}^n = 0$,  the $h$-order terms in the
parenthesis vanish and taking in the summation over $i$ and
considering periodicity of the Volterra lattice one can easily see
that the second terms of order $h$ cancel and the quadratic Poisson
brackets (11) are preserved by the symplectic Euler method.

Because the second order Lobatto IIIA-B method is a composition of
symplectic Euler methods with  step sizes $h/2$, it preserves the
quadratic Poisson bracket of the Volterra lattice. But higher order
Lobatto IIIA-B methods can not preserve the Poisson structure,
because they can not be  written as combination of symplectic Euler
method. Only diagonally implicit partitioned Runge-Kutta methods can
be written as combination of symplectic Euler method (see pp. 180,
\cite{hairer02gi}).

Recently finite dimensional systems with a Poisson structure arising
after semi-discretization of certain partial differential equations
have been integrated by splitting methods. An example of this is the
preservation of the Lie-Poisson structure of Landau-Lifschitz
lattice in partitioned form using a staggered scheme which
corresponds to the second order Lobatto IIIA-B \cite{frank97gif}. It
was shown in \cite{faou04pig} that the finite dimensional system
which arises using a variational approximation of the time-dependent
Schr\"odinger equation by Gaussian wave packets inherits a Poisson
structure and various splitting methods were considered for its
integration.

The linear integrals are preserved exactly by all Runge-Kutta
methods. The implicit mid-point rule preserves the quadratic
integrals exactly assuming that the underlying system of nonlinear
equations is solved within the machine accuracy. The Lobatto IIIA-B
methods preserve only the quadratic integrals of the form $Q=p^TDq$,
where $D$ is an arbitrary matrix of appropriate dimension  (Ch. 4,
\cite{hairer02gi}). Unfortunately the quadratic integral $I_q$ (13)
of the splitted Volterra lattice is not in this form. Higher order
polynomial integrals like the cubic integral $I_c$ and
non-polynomial integrals like the Casimir function $H_0$ are not
preserved exactly by the symplectic Euler and Lobatto IIIA-B
methods.

For the periodic Volterra lattice we have used the following initial
condition
$$
y(x_i) = 1 + \frac{1}{2m^2}\sec\!{\rm h} ^2 (x_i), \quad x_i= -1+
(i-1)\frac{1}{2m}, \quad i=1,\ldots, m.
$$
All computations are done with a constant time step $\Delta t = 0.1$
over the time interval $t\in[0,2000]$ for a Volterra lattice of
dimension $m=20,40,80$. The errors in the Hamiltonians and conserved
quantities are given in Table 1 in the mean square root norm
$\sqrt{\sum_{i=1}^N (I^i-I_0)^2}/N$, where $I^i$ denote the computed
Hamiltonians or first integrals at time step $t_i$ and $N$ is the
number of time steps.

For all values of $m$ and $\Delta t$, the Hamiltonian $H_1$ is
preserved with almost the same high accuracy for both methods. The
Casimir $H_0$, the quadratic and cubic first integrals can not be
preserved exactly, but the errors do not grow with time as in
non-symplectic methods.  Similar numerical results are obtained for
the Toda lattice (see pp.
 385-386 \cite{hairer02gi} and pp. 430-431  \cite{hairer03gni}). It
was shown in \cite{hairer02gi}, Theorem 3.1, pp. 353, that for
completely integrable systems, the symplectic integrators preserve
the first integrals  over long-time  with an error ${\cal O}(\Delta
t^p)$, where $p$ denotes the order of the method. The Lobatto IIIA-B
method results in smaller errors than the symplectic Euler method
because it is a second order accurate method. One can also observe
that the Casimir $H_0$ is preserved slightly better than the
quadratic and cubic first integrals

%The different form of the fluctuations of the errors of $H_0$, the
%quadratic and cubic integrals over time in Figures 1 and 2 for the
%symplectic Euler and Lobatto IIIA-B method are due to the difference
%of the modified equations for both methods, because the Hamiltonians
%and first integrals also satisfy similar modified equations like the
%Volterra lattice equation.
\newpage
\begin{center}
\begin{figure}[ht]
\epsfxsize=3.5in \centerline{\epsfbox{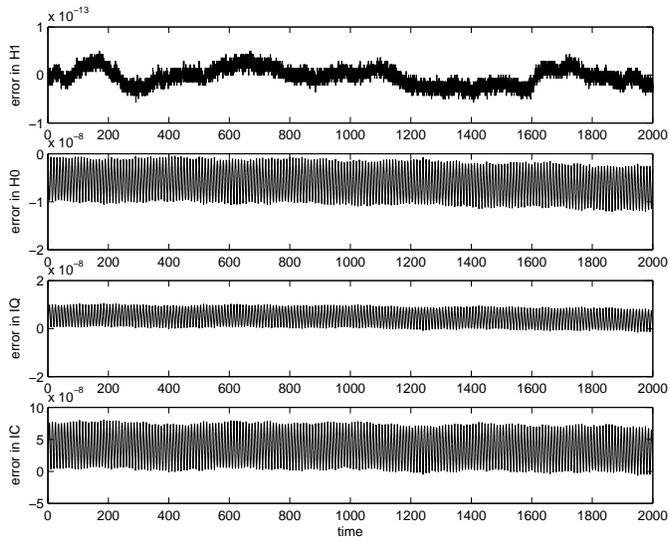}} \caption{Errors
in the Hamiltonians and first integrals: symplectic Euler method,
$m=40$} \label{se}
\end{figure}
\end{center}

\begin{center}
\begin{figure}[ht]
\epsfxsize=3.5in \centerline{\epsfbox{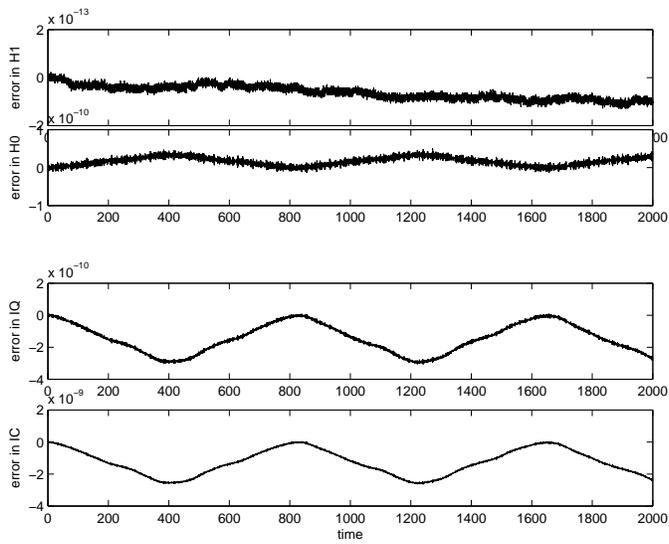}} \caption{Errors
in the  Hamiltonians and first integrals:
 Lobatto IIIA-B method, $m=40$} \label{lo}
\end{figure}
\end{center}
\newpage
\begin{table}[htb]
\caption{Average errors of the Hamiltonians and first
integrals}
\begin{center}
\begin{tabular}{|l|l|cccc|cccc|}\hline
 & &  \multicolumn{4}{c|}{symplectic Euler method} & \multicolumn{4}{c|}
{Lobatto
 IIIA-B method}\\ \hline
 $m$  & $ \Delta t$ & $H_1$ & $H_0$ & $I_q$ & $ I_c$
& $H_1$ & $H_0$ & $I_q$ & $ I_c$\\ \hline
20 & 0.2 & 2.589 -16 & 7.021 -08 & 6.200 -08 & 2.380 -07 & 3.042 -16   & 2.698 -11   & 1.879 -10   & 1.597 -09\\
20 & 0.1 & 2.740 -16 & 1.240 -08 & 1.072 -08 & 4.100 -08 & 2.526
-16 &   5.157 -12 &   3.268 -11 &   2.765 -10  \\
 20 &
0.05 & 1.652 -16   & 2.229 -09 &   1.853 -09 & 7.036 -09 &
1.609 -16 &  9.373 -13 &   5.793 -12 &   4.896 -11 \\
\hline
40 & 0.2 & 3.016 -16 &  2.608 -10 &   2.158 -10 &   1.746 -09 &   1.752  -16 &  7.557 -13 &  6.870 -12 & 6.036 -11\\
40 & 0.1 & 1.283 -16 &  4.820 -11 &  4.021 -11 &  3.218 -10 &  4.824 -16 &  1.401 -13 &   1.202 -12 &   1.053  -11\\
40 & 0.05 & 2.309 -16 &  9.300 -12 &  7.863 -12 &  6.163 -11 & 2.562 -16 &  2.540 -14 &  2.115 -13 &  1.857 -12\\
\hline
80 & 0.2 & 2.903 -16 &  7.256 -12 &   7.162 -12 &   5.650 -11 & 9.837 -16 &  1.409 -14 &  1.396 -13 &  1.229 -12 \\
80 & 0.1 & 5.237 -16 &   1.314 -12 &  1.296 -12 &   1.017 -11 & 2.244 -16 & 2.968 -15  & 2.342 -14 &  2.137 -13\\
80 & 0.05 & 6.117 -16 &   2.439 -13  &   2.396 -13 &   1.864 -12 & 1.547 -16 &  6.175 -16 & 3.830 -15 &  3.745 -14 \\
\hline
\end{tabular}
\end{center}
\end{table}

%Near-preservation of first integrals. All first integrals are
%preserved with an accuracy of order $p$. For symplectic Eueler
%method $p=1$, for the second order Lobatto method $p=2$. Compute
%with the step size h=0.05. There is no drift the error for first
%integrals in contrast to standard integrators like Runge-Kutta
%methods, see for the Toda lattice Hairer pp. 385-386 and Hairer,
%pp 353, Theorem 3.1. Symplectic methods preserve the Hamiltonians
%and first integrals over exponentially long time interval  (? ).
%We have
%$$
%H(y_n) = H(y_0) + {\cal O\/}(h^p)
%$$
%for symplectic integrators, where $m$ denotes the order of the
%methods. Non-symplectic methods show a drift in the Hamiltonians
%and first integrals in long time integration unless very small
%sizes are used.
%$$
%H(y_n) = H(y_0) + {\cal O\/}(nh^p)
%$$

\section{Conclusion}

We have shown that the symplectic Euler method preserves the
quadratic Poisson structure of the periodic Volterra lattice.  The
numerical results show excellent long time preservation of the
Hamiltonian, Casimirs and the first integrals. Because of the
singularity of the structure matrices we can obtain only local
results for the backward error analysis in contrast to the global
results obtained for Poisson systems with invertible structure
matrices like in (Sec. IX.3.3, pp.297 \cite{hairer02gi}).

\section{Acknowledgements}

The second author  acknowledges the support of Swiss National
Science Foundation and is grateful to Ernst Hairer and Gerhard
Wanner for their hospitality during his stay at Universit\'{e} de
Gen\`{e}ve. The authors thank to the referees for helpful comments
and suggestions.

\bibliographystyle{plain}

\end{document}